\documentclass{article}
\usepackage[utf8]{inputenc}

\title{Properties of Hamiltonian Circuits in Rectangular Grids}
\author{R\"udiger Jehn }
\date{March 2021}

\usepackage{tikz}
\usepackage{manfnt}
\usepackage{diagbox}
\usepackage{color, colortbl}
\usepackage[colorlinks,citecolor=blue,urlcolor=blue,bookmarks=false,hypertexnames=true]{hyperref} 
\definecolor{Gray}{gray}{0.9}

\begin{document}

\maketitle

{\bf Abstract}\\

We present properties and invariants of Hamiltonian circuits in rectangular grids. It is proved that all circuits on a $2n \times 2n$ chessboard have at least $4n$ turns and at least $2n$ straights if $n$ is even and $2n+2$ straights if $n$ is odd. The minimum number of turns and straights are presented and proved for circuits on an $n \times (n+1)$ chessboard. For the general case of an $n \times m$ chessboard similar results are stated but not all proofs are given.

\section{Introduction}

\begin{figure}[h]
\begin{center}
\begin{tikzpicture}

\draw (-2,-2) -- (2,-2);
\draw (-2,-1) -- (2,-1);
\draw (-2,0) -- (2,0);
\draw (-2,1) -- (2,1);
\draw (-2,2) -- (2,2);

\draw (-2,-2) -- (-2,2);
\draw (-1,-2) -- (-1,2);
\draw (0,-2) -- (0,2);
\draw (1,-2) -- (1,2);
\draw (2,-2) -- (2,2);

\draw[blue, thick] (-1.5,-1) -- (-1.5,1);
\draw[blue, thick] (-1.5,1) .. controls +(up:3mm) and +(left:3mm) .. (-1, 1.5);
\draw[blue, thick] (-1, 1.5) -- (1,1.5);
\draw[blue, thick] (1,1.5) .. controls +(right:3mm) and +(up:3mm) .. (1.5, 1);
\draw[blue, thick] (1.5,1) .. controls +(down:3mm) and +(right:3mm) .. (1, 0.5);
\draw[blue, thick] (1, 0.5) -- (0,0.5);
\draw[blue, thick] (0,0.5) .. controls +(left:3mm) and +(up:3mm) .. (-0.5, 0);
\draw[blue, thick] (-0.5,0) .. controls +(down:3mm) and +(left:3mm) .. (0, -0.5);
\draw[blue, thick] (0, -0.5) -- (1, -0.5);
\draw[blue, thick] (1,-0.5) .. controls +(right:3mm) and +(up:3mm) .. (1.5, -1);
\draw[blue, thick] (1.5,-1) .. controls +(down:3mm) and +(right:3mm) .. (1, -1.5);
\draw[blue, thick] (-1, -1.5) -- (1, -1.5);
\draw[blue, thick] (-1, -1.5) .. controls +(left:3mm) and +(down:3mm) .. (-1.5,-1);

\end{tikzpicture}
\caption{A Hamiltonian (or rook) circuit on a $4 \times 4$ chessboard.}
\label{fig1}
\end{center}
\end{figure}
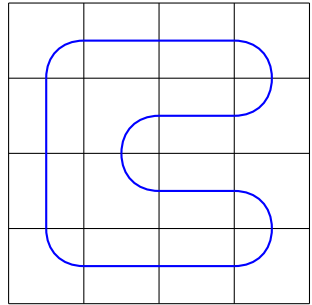

A Hamiltonian circuit in a rectangular grid, often also called a {\it rook circuit} or Wazir tour \cite{jeliss} is a circuit on an $n \times m$ chessboard that passes through all squares without crossing. Sometimes the term Hamiltonian circuit is also used for a circuit on a chessboard of a knight \cite{shufelt}, but in this paper we limit ourselves to moves that a rook is allowed to make. A complete circuit is obviously only possible if either $n$ or $m$ or both are even numbers. Figure \ref{fig1} shows one possible circuit on a $4 \times 4$ chessboard.

\vspace{-5mm}
There are many ways to draw a Hamiltonian circuit. Stoyan and Strehly~\cite{stoyan} and four years later Kloczkowski \cite{kloczkowski} developed algorithms to enumerate how many Hamiltonian circuits there are for any $n \times m$ grid. The main focus of this paper is to find the circuits with the smallest or largest number of turns. In the example above there are 8 squares where the driver needs to turn. The first question is: How many turns must a rook circuit at least have? \\

The next question is: What is the maximum number of turns in a rook circuit or equivalently what is the minimum number of straights? We start with proving some basic properties of any Hamiltonian circuit in a rectangular grid.

\section{Number of corner points inside a Hamiltonian circuit}

We call the centres of each square of the chessboard grid points and the corners of the squares are called corner points. \\

\subsection{Lemma 1}
The number of corner points inside a Hamiltonian circuit on $n \times m$ grid is 
$$\frac{n \times m}{2} - 1.$$

Proof:\\

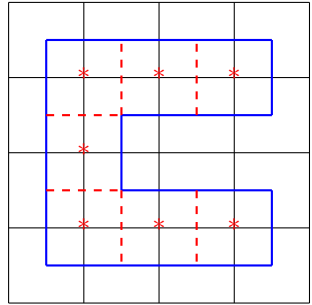
\begin{figure}[h]
\begin{center}
\begin{tikzpicture}

\draw (-2,-2) -- (2,-2);
\draw (-2,-1) -- (2,-1);
\draw (-2,0) -- (2,0);
\draw (-2,1) -- (2,1);
\draw (-2,2) -- (2,2);

\draw (-2,-2) -- (-2,2);
\draw (-1,-2) -- (-1,2);
\draw (0,-2) -- (0,2);
\draw (1,-2) -- (1,2);
\draw (2,-2) -- (2,2);

\draw[blue, thick] (-1.5,-1.5) -- (-1.5,1.5);
\draw[blue, thick] (-1.5, -1.5) -- (1.5,-1.5);
\draw[blue, thick] (-1.5, 1.5) -- (1.5,1.5);
\draw[blue, thick] (-0.5, -0.5) -- (1.5,-0.5);
\draw[blue, thick] (-0.5, 0.5) -- (1.5,0.5);
\draw[blue, thick] (1.5, 1.5) -- (1.5,0.5);
\draw[blue, thick] (1.5, -0.5) -- (1.5,-1.5);
\draw[blue, thick] (-0.5, -0.5) -- (-0.5,0.5);

\draw[red, dashed, thick] (-1.5, 0.5) -- (-0.5, 0.5);
\draw[red, dashed, thick] (-1.5, -0.5) -- (-0.5, -0.5);
\draw[red, dashed, thick] (-0.5, 0.5) -- (-0.5, 1.5);
\draw[red, dashed, thick] (0.5, 0.5) -- (0.5, 1.5);
\draw[red, dashed, thick] (-0.5, -0.5) -- (-0.5, -1.5);
\draw[red, dashed, thick] (0.5, -0.5) -- (0.5, -1.5);
\draw[red] (-1,1) -- (-1,1)   node[anchor=center]  {*};
\draw[red] (0,1) -- (0,1)   node[anchor=center]  {*};
\draw[red] (1,1) -- (1,1)   node[anchor=center]  {*};
\draw[red] (-1,0) -- (-1,0)   node[anchor=center]  {*};
\draw[red] (-1,-1) -- (-1,-1)   node[anchor=center]  {*};
\draw[red] (0,-1) -- (0,-1)   node[anchor=center]  {*};
\draw[red] (1,-1) -- (1,-1)   node[anchor=center]  {*};

\end{tikzpicture}
\caption{A Hamiltonian (or rook) circuit on a 4 x 4 chessboard has always 7 corner points inside the circuit.}
\label{fig:pick}
\end{center}
\end{figure}

Lemma 1 follows directly from Pick's Theorem \cite{pick} which states that the area inside a polygon that is defined by a sequence of grid points is the number of interior points plus half the number of points on the polygon minus 1. If we apply it to any Hamiltonian circuit, the number of interior points is zero and the number of points on the boundary is n $\times$ $m$ because the polygon passes through all grid points and hence the area inside the circuit is $\frac{n \times m}{2} - 1$. But since the chessboard can be split into unit squares around the corner points (see red dashed squares in Figure \ref{fig:pick}), the number of corner points inside the polygon is equal to the area of the polygon.\\

q.e.d.

\section{Invariants of Hamiltonian circuits}

Each turn in a square of a chessboard is either clockwise or counter-clockwise. We will define, without loss of generality, the direction of the circuit such that the turn in the upper left corner will be counter-clockwise. This is in agreement with the way a 400-m track is run in a stadium. However, it is opposite to most of the major European formula 1 tracks which are mostly clockwise.\\

\subsection{Lemma 2}

Lemma 2: Let $c$ be the number of squares of a chessboard where the rook turns counter-clockwise and $d$ the number of squares where the rook turns clockwise, then $c - d = 4$.\\

This lemma is pretty straightforward to prove. Since the path is not allowed to cross, the sum of all turns (either $+90^\circ$ or $-90^\circ$) of the final circuit must add up to $360^\circ$. Hence we need 4 more turns in counter-clockwise direction than in clockwise direction. A simple observation is the following: Since the turn in the upper left corner is defined to be counter-clockwise then right of the path is the "outside" and left of the path is the "inside". In order to have the edges of the chessboards on the outside of the circuit, there must be turns in counter-clockwise direction in all 4 corners.

\subsection{Lemma 3}

Lemma 3: A rook circuit covering a complete chessboard has an even number of turns in every row or column.\\

We will prove the lemma for rows, as the proof for columns is completely identical. Let us take the $k$ upper rows of an $n \times m$ chessboard. The number of turns in this $k \times m$ field must be even, because on each loop on which the rook enters and exits the field the rook changes $x$ times from vertical to horizontal motion and also $x$ times from horizontal to vertical motion. This means the rook does $2x$ turns. It can enter and leave the $k \times m$ field various times, but the sum of its turns will always be an even number.\\

And if the number of turns in an $k \times m$ field is even and it is also even in an $(k+1) \times m$ field, it must be even in row $k+1$.\\

q.e.d.\\

Corollary: It follows directly from Lemma 3 that the number of straight crossings of a square is even (odd) in a row or column if the number of squares in the row or column is even (odd).\\

\subsection{Lemma 4}

Lemma 4: Let $s_k \in \{0, 1\}^{2l}$ denote the sequence of entries (1) and exits (0) from row $k$ to the next row $k+1$ as enumerated from left to right then $s_k$ will always be of the form $(10)^l$.\\

This means that the first square from the left in row $k+1$ that is connected with the row above, has an entry from above. And the second square from the left in row $k+1$ that is connected with the row above, has an exit to above. And this sequence keeps repeating from left to right.\\

Proof:\\

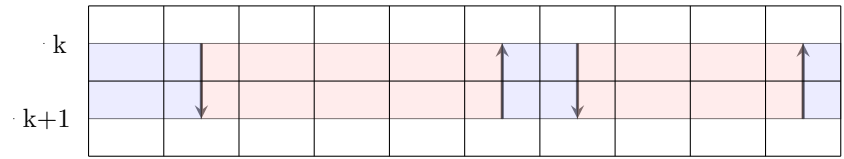
\begin{figure}[h]
\begin{center}
\begin{tikzpicture}

\draw (-0.6,1.5) -- (-0.6,1.5)   node[anchor=west]  {k};
\draw (-1,0.5) -- (-1,0.5)   node[anchor=west]  {k+1};

\draw[-stealth, very thick] (1.5, 1.5) -- (1.5, 0.5);
\draw[-stealth, very thick] (6.5, 1.5) -- (6.5, 0.5);
\draw[-stealth, very thick] (5.5, 0.5) -- (5.5, 1.5);
\draw[-stealth, very thick] (9.5, 0.5) -- (9.5, 1.5);
\filldraw[fill=blue!15, opacity=0.5] (0,0.5) rectangle (1.5,1.5);
\filldraw[fill=red!15, opacity=0.5] (1.5,0.5) rectangle (5.5,1.5);
\filldraw[fill=blue!15, opacity=0.5] (5.5,0.5) rectangle (6.5,1.5);
\filldraw[fill=red!15, opacity=0.5] (6.5,0.5) rectangle (9.5,1.5);
\filldraw[fill=blue!15, opacity=0.5] (9.5,0.5) rectangle (10,1.5);


\draw (0,0) -- (10,0);
\draw (0,1) -- (10,1);
\draw (0,2) -- (10,2);
\draw (0,0) -- (0,2);
\draw (1,0) -- (1,2);
\draw (2,0) -- (2,2);
\draw (3,0) -- (3,2);
\draw (4,0) -- (4,2);
\draw (5,0) -- (5,2);
\draw (6,0) -- (6,2);
\draw (7,0) -- (7,2);
\draw (8,0) -- (8,2);
\draw (9,0) -- (9,2);
\draw (10,0) -- (10,2);

\end{tikzpicture}
\caption{The blue shaded area shows the outside of a circuit and the red shaded area the inside. Inside and outside must take turns along the boundary between row $k$ and row $k+1$ and hence also entries and exits must takes turns.}
\label{fig:1010}
\end{center}
\end{figure}

The first observation is that the number of entries between row $k$ and $k+1$ must be the same as the number of exits in order to have a closed circuit. Hence the length of $s_k$ is even. Due to the definition of the counter-clockwise orientation of the circuit, the right side of the path is the "outside" and the left side of the path is the "inside". The first connection (starting from the left) between row $k$ and row $k+1$ must be an entry because in this way, the left border of the chessboard is outside the circuit. Figure \ref{fig:1010} illustrates this geometry by colouring the outside of the circuit in blue and the inside in red.\\

Following the same argument, the second connection between row $k$ and row $k+1$ (drawn in column 6 in Figure \ref{fig:1010}) must be an exit, because the inside of the circuit is to the left of the connection. Inside and outside of the circuit must take turns along the boundary between row $k$ and row $k+1$ and hence also entries and exits must takes turns, and hence $s_k$ must have the form $(10)^l$.\\

q.e.d.\\

(Note: a symmetrical statement holds obviously also for any column of the chessboard.)

\subsection{Lemma 5}

Lemma 5: Let $c$ be the number of squares of a row of a chessboard with $m$ columns where the rook turns counter-clockwise, $d$ the number of squares in this row where the rook turns clockwise, $e$ the number of squares where the rook moves horizontally and $f$ the number of squares where the rook moves vertically, then $(c + e - d - f) \bmod 4 = m \bmod 4$.\\

Proof:\\

As often the rook enters row $k$ it has to leave the row again. We consider the 4 possible cases:

\begin{enumerate}
    \item The rook enters row $k$ from row $k$-1 and leaves towards row $k$-1
    \item The rook enters row $k$ from row $k$-1 and leaves towards row $k$+1
    \item The rook enters row $k$ from row $k$+1 and leaves towards row $k$-1
    \item The rook enters row $k$ from row $k$+1 and leaves towards row $k$+1
\end{enumerate}

It is easy to verify that in case 2 and 3 where the rook traverses row $k$ (but possibly moving a few squares horizontally in that row), the number of squares in this row passage where the rook turns clockwise plus the number of vertical cell crossings is exactly 1. In case 1 and 4 this sum is either 0 or 2. But case 2 (passage from above to below) must happen as many times as case 3 (passage from below to above). Hence the sum of the number of squares of all row passages (cases 2 and 3 together) where the rook turns clockwise and the number of vertical cell crossings is an even number. Adding the contributions from cases 1 and 4 will not change the parity of this sum and hence the sum of $d$ and $f$ is even.\\

It is a basic rule if $a + b = m$ and $b$ is even then $(a - b) \bmod 4 = m \bmod 4$. Hence we can conclude that $(c + e - d - f) \bmod 4 = m \bmod 4$.\\

q.e.d.\\

\section{Minimum number of turns in a rook circuit on a $2n \times 2n$ chessboard}

We will prove that every rook circuit in a $2n \times 2n$ chessboard has at least 4n turns.

\subsection{Lemma 6}

Lemma 6: Each $m \times m$ corner of a chessboard has at least $m$ turns.\\ 

This lemma will be proved with induction. We initialise the induction proof with $m=1$. There is always a turn in each corner of the chessboard and therefore the claim holds for $m=1$.\\

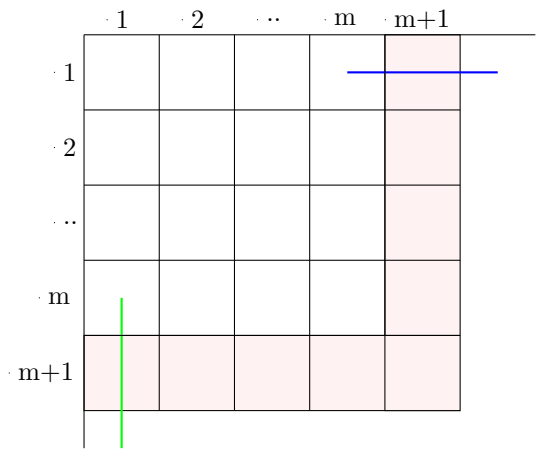
\begin{figure}[h]
\begin{center}
\begin{tikzpicture}

\draw (-1.7,3.2) -- (-1.7,3.2)   node[anchor=west]  {1};
\draw (-0.7,3.2) -- (-0.7,3.2)   node[anchor=west]  {2};
\draw (0.3,3.2) -- (0.3,3.2)   node[anchor=west]  {..};
\draw (1.2,3.2) -- (1.2,3.2)   node[anchor=west]  {m};
\draw (2.0,3.2) -- (2.0,3.2)   node[anchor=west]  {m+1};
\draw (-2.4,2.5) -- (-2.4,2.5)   node[anchor=west]  {1};
\draw (-2.4,1.5) -- (-2.4,1.5)   node[anchor=west]  {2};
\draw (-2.4,0.5) -- (-2.4,0.5)   node[anchor=west]  {..};
\draw (-2.6,-0.5) -- (-2.6,-0.5)   node[anchor=west]  {m};
\draw (-3.0,-1.5) -- (-3.0,-1.5)   node[anchor=west]  {m+1};

\filldraw[fill=red!5] (2,-2) rectangle (3,3);
\filldraw[fill=red!5] (-2,-2) rectangle (3,-1);

\draw[blue, thick] (1.5, 2.5) -- (3.5, 2.5);
\draw[green, thick] (-1.5, -2.5) -- (-1.5, -0.5);

\draw (-2,-1) -- (3,-1);
\draw (-2,0) -- (3,0);
\draw (-2,1) -- (3,1);
\draw (-2,2) -- (3,2);
\draw (-2,3) -- (4,3);

\draw (-2,-2.5) -- (-2,3);
\draw (-1,-2) -- (-1,3);
\draw (0,-2) -- (0,3);
\draw (1,-2) -- (1,3);
\draw (2,-2) -- (2,3);

\end{tikzpicture}
\caption{Left upper corner of the chessboard extended by one row and one column}
\label{fig2}
\end{center}
\end{figure}

Now we assume that the lemma holds for $m$, i.e.\ we have at least $m$ turns in each $m \times m$ corner. And let us assume that the $(m+1) \times (m+1)$ part of the chessboard will have no turns in row $m+1$ and in column $m+1$ (the red area in the Figure \ref{fig2}). Without loss of generality we examine the left upper corner of the chessboard. Since we have no turns in column $m+1$, we need to have a horizontal connection between squares $a_{1,m}$ and $a_{1,m+2}$ (drawn in blue). Similarly we need to have a vertical connection between squares $a_{m,1}$ and $a_{m+2,1}$ (drawn in green).

\vspace{-2mm}

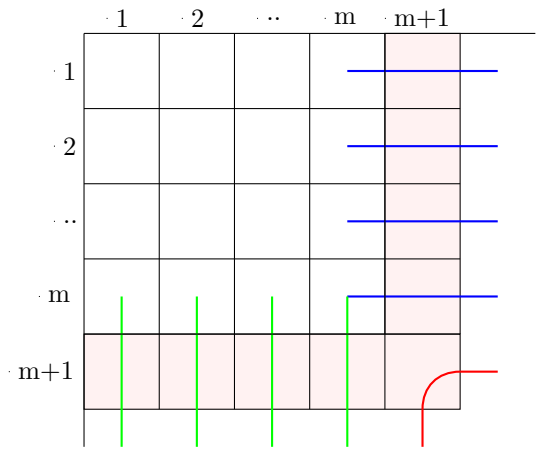
\begin{figure}
\begin{center}
\begin{tikzpicture}

\draw (-1.7,3.2) -- (-1.7,3.2)   node[anchor=west]  {1};
\draw (-0.7,3.2) -- (-0.7,3.2)   node[anchor=west]  {2};
\draw (0.3,3.2) -- (0.3,3.2)   node[anchor=west]  {..};
\draw (1.2,3.2) -- (1.2,3.2)   node[anchor=west]  {m};
\draw (2.0,3.2) -- (2.0,3.2)   node[anchor=west]  {m+1};
\draw (-2.4,2.5) -- (-2.4,2.5)   node[anchor=west]  {1};
\draw (-2.4,1.5) -- (-2.4,1.5)   node[anchor=west]  {2};
\draw (-2.4,0.5) -- (-2.4,0.5)   node[anchor=west]  {..};
\draw (-2.6,-0.5) -- (-2.6,-0.5)   node[anchor=west]  {m};
\draw (-3.0,-1.5) -- (-3.0,-1.5)   node[anchor=west]  {m+1};

\filldraw[fill=red!5] (2,-2) rectangle (3,3);
\filldraw[fill=red!5] (-2,-2) rectangle (3,-1);

\draw[blue, thick] (1.5, 2.5) -- (3.5, 2.5);
\draw[blue, thick] (1.5, 1.5) -- (3.5, 1.5);
\draw[blue, thick] (1.5, 0.5) -- (3.5, 0.5);
\draw[blue, thick] (1.5, -0.5) -- (3.5, -0.5);
\draw[green, thick] (-1.5, -2.5) -- (-1.5, -0.5);
\draw[green, thick] (-0.5, -2.5) -- (-0.5, -0.5);
\draw[green, thick] (0.5, -2.5) -- (0.5, -0.5);
\draw[green, thick] (1.5, -2.5) -- (1.5, -0.5);
\draw[red, thick] (2.5,-2.5) -- (2.5,-2.0);
\draw[red, thick] (2.5,-2.0) .. controls +(up:3mm) and +(left:3mm) .. (3.0, -1.5);
\draw[red, thick] (3.0, -1.5) -- (3.5, -1.5);

\draw (-2,-1) -- (3,-1);
\draw (-2,0) -- (3,0);
\draw (-2,1) -- (3,1);
\draw (-2,2) -- (3,2);
\draw (-2,3) -- (4,3);

\draw (-2,-2.5) -- (-2,3);
\draw (-1,-2) -- (-1,3);
\draw (0,-2) -- (0,3);
\draw (1,-2) -- (1,3);
\draw (2,-2) -- (2,3);

\end{tikzpicture}
\caption{There is a dilemma if you want to fill the red squares with straight lines only.}
\label{fig3}
\end{center}
\end{figure}

But this infers that also the squares $a_{2,m}$ and $a_{2,m+2}$ must be connected by a straight line. And this propagates all the way down to $a_{m,m}$ and $a_{m,m+2}$. The same argument holds for the squares $a_{m,2}$ and $a_{m+2,2}$ all the way to $a_{m,m}$ and $a_{m+2,m}$. So we end up with square $a_{m+1,m+1}$ that cannot connect to $a_{m,m+1}$ nor to $a_{m+1,m}$ (see Figure \ref{fig3}).\\


The only possibility for square $a_{m+1,m+1}$ is a to have a turn (like drawn in red). This is a contradiction to the assumption that there are no turns in column and row $m+1$. But if we have a turn in the red area, and we have at least $n$ turns in the $m \times m$ quadrant, we will have at least $m+1$ turns in a $(m+1) \times (m+1)$ quadrant which completes our induction proof. \\

\subsection{Final synthesis}

Now we divide the $2n \times 2n$ chessboard into 4 quadrants with size $n \times n$. Since each quadrant has at least $n$ turns, the minimum number of turns in all four quadrants is $4n$.\\

q.e.d.


\section{Minimum number of straights in a rook circuit on a $2n \times 2n$ chessboard}

In the OEIS sequence A085622 \cite{A085622} a formula is given for the maximal number of segments (equivalently, corners) in a rook circuit of a $2n \times 2n$ board. If we subtract this number from the total number of squares, i.e.\ $4n^2$ we get this formula for the minimum number of straights:

$$
a(n) = \left\{
    \begin{array}{ll}
        2n & \mbox{if n is even} \\
        2n + 2 & \mbox{if n is odd and} > 1 
    \end{array}
\right.
$$

This will be proved here.

\subsection{Lemma 7}

Lemma 7: Each $2m \times 2m$ corner of a chessboard has at least $m$ straights.\\ 

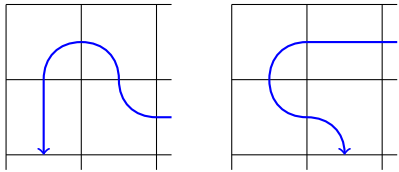
\begin{figure}[h]
\begin{center}
\begin{tikzpicture}

\draw (-2,0) -- (0.2,0);
\draw (-2,1) -- (0.2,1);
\draw (-2,2) -- (0.2,2);

\draw (-2,-0.2) -- (-2,2);
\draw (-1,-0.2) -- (-1,2);
\draw (0,-0.2) -- (0,2);

\draw[blue, thick] (0.2,0.5) -- (0,0.5);
\draw[blue, thick] (0,0.5) .. controls +(left:3mm) and +(down:3mm) .. (-0.5, 1);
\draw[blue, thick] (-0.5, 1) .. controls +(up:3mm) and +(right:3mm) .. (-1, 1.5);
\draw[blue, thick] (-1, 1.5) .. controls +(left:3mm) and +(up:3mm) .. (-1.5, 1);
\draw[->, blue, thick] (-1.5, 1) -- (-1.5,0);

\draw (1,0) -- (3.2,0);
\draw (1,1) -- (3.2,1);
\draw (1,2) -- (3.2,2);

\draw (1,-0.2) -- (1,2);
\draw (2,-0.2) -- (2,2);
\draw (3,-0.2) -- (3,2);

\draw[blue, thick] (3.2, 1.5) -- (2, 1.5);
\draw[blue, thick] (2, 1.5) .. controls +(left:3mm) and +(up:3mm) .. (1.5, 1);
\draw[blue, thick] (1.5, 1) .. controls +(down:3mm) and +(left:3mm) .. (2, 0.5);
\draw[->, blue, thick] (2, 0.5) .. controls +(right:3mm) and +(up:3mm) .. (2.5, 0);

\end{tikzpicture}
\caption{The rook either leaves (left) or enters (right) the upper left corner on a straight line.}
\label{fig2by2}
\end{center}
\end{figure}

This lemma will also be proved with induction very similarly to the proof of lemma 6. We initialise the induction proof with $m=1$. Figure \ref{fig2by2} shows the two basic possibilities with one straight line in a 2 x 2 corner. There is also the possibility to have two straight lines which is not plotted here. The only solution with no straight lines would be a circle, which cannot be connected with a complete rook circuit.

\vspace{-5mm}
Now we assume that we have a $2m \times 2m$ corner containing $m$ straight lines. And let us assume that the $2(m+1) \times 2(m+1)$ extended part of the chessboard will have no straight lines in the rows and columns $m+1$ and $m+2$.\\ 

We start investigating square $a_{1,m+1}$. Since we cannot have a straight line in it, the only path could be case 1 \manquartercircle \ or case 2 \manrotatedquartercircle. Lets assume the first case and look what that means for the square square $a_{2,m+1}$. In this square underneath we can only have either case 1.1 \rotatebox[origin=c]{180}{\manrotatedquartercircle} \ or case 1.2 \rotatebox[origin=c]{180}{\manquartercircle} if we want to avoid a straight line. In case 1.1 this means we will automatically have \manrotatedquartercircle \ in square $a_{1,m+2}$ and \rotatebox[origin=c]{180}{\manquartercircle} in square $a_{2,m+2}$ (see Figure \ref{case11}). In case 1.2 we will automatically have \rotatebox[origin=c]{180}{\manrotatedquartercircle} in square $a_{2,m+2}$ and \manrotatedquartercircle \ in square $a_{1,m+2}$ and (see blue path in Figure \ref{case2}). \\

In case 2 we will automatically have a \manquartercircle \ in square $a_{1,m+2}$ to avoid the straight line and this forces a \rotatebox[origin=c]{180}{\manquartercircle} in square $a_{2,m+2}$. And the only choice for square $a_{2,m+1}$ will be \rotatebox[origin=c]{180}{\manrotatedquartercircle} (see red path in Figure \ref{case2}). \\

\vspace{-2mm}

\begin{figure}[h]
\begin{center}
\begin{tikzpicture}

\draw (-1.7,3.2) -- (-1.7,3.2)   node[anchor=west]  {1};
\draw (-0.7,3.2) -- (-0.7,3.2)   node[anchor=west]  {2};
\draw (0.3,3.2) -- (0.3,3.2)   node[anchor=west]  {..};
\draw (1.2,3.2) -- (1.2,3.2)   node[anchor=west]  {m};
\draw (2.0,3.2) -- (2.0,3.2)   node[anchor=west]  {m+1};
\draw (3.0,3.2) -- (3.0,3.2)   node[anchor=west]  {m+2};
\draw (-2.4,2.5) -- (-2.4,2.5)   node[anchor=west]  {1};
\draw (-2.4,1.5) -- (-2.4,1.5)   node[anchor=west]  {2};

\filldraw[fill=red!5] (2,1) rectangle (4,3);

\draw (-2,1) -- (4,1);
\draw (-2,2) -- (4,2);
\draw (-2,3) -- (4.5,3);

\draw (-2,0.8) -- (-2,3);
\draw (-1,1) -- (-1,3);
\draw (0,1) -- (0,3);
\draw (1,1) -- (1,3);
\draw (2,1) -- (2,3);
\draw (3,1) -- (3,3);

\draw[blue, thick] (1.5, 2.5) -- (2, 2.5);
\draw[blue, thick] (2, 2.5) .. controls +(right:3mm) and +(up:3mm) .. (2.5, 2);
\draw[blue, thick] (2.5, 2) .. controls +(down:3mm) and +(right:3mm) .. (2, 1.5);
\draw[blue, thick] (2, 1.5) -- (1.5, 1.5);

\draw[blue, thick] (4.5, 2.5) -- (4, 2.5);
\draw[blue, thick] (4, 2.5) .. controls +(left:3mm) and +(up:3mm) .. (3.5, 2);
\draw[blue, thick] (3.5, 2) .. controls +(down:3mm) and +(left:3mm) .. (4, 1.5);
\draw[blue, thick] (4, 1.5) -- (4.5, 1.5);

\end{tikzpicture}
\vspace{-1mm}
\caption{Case 1.1}
\label{case11}
\end{center}
\end{figure}
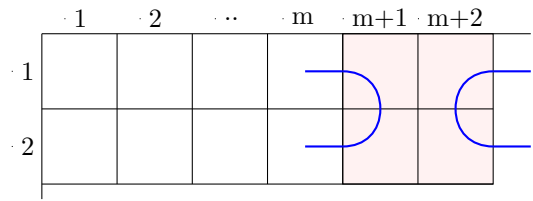

\begin{figure}[h]
\begin{center}
\begin{tikzpicture}

\draw (-1.7,3.2) -- (-1.7,3.2)   node[anchor=west]  {1};
\draw (-0.7,3.2) -- (-0.7,3.2)   node[anchor=west]  {2};
\draw (0.3,3.2) -- (0.3,3.2)   node[anchor=west]  {..};
\draw (1.2,3.2) -- (1.2,3.2)   node[anchor=west]  {m};
\draw (2.0,3.2) -- (2.0,3.2)   node[anchor=west]  {m+1};
\draw (3.0,3.2) -- (3.0,3.2)   node[anchor=west]  {m+2};
\draw (-2.4,2.5) -- (-2.4,2.5)   node[anchor=west]  {1};
\draw (-2.4,1.5) -- (-2.4,1.5)   node[anchor=west]  {2};

\filldraw[fill=red!5] (2,1) rectangle (4,3);

\draw (-2,1) -- (4,1);
\draw (-2,2) -- (4,2);
\draw (-2,3) -- (4.5,3);

\draw (-2,0.8) -- (-2,3);
\draw (-1,1) -- (-1,3);
\draw (0,1) -- (0,3);
\draw (1,1) -- (1,3);
\draw (2,1) -- (2,3);
\draw (3,1) -- (3,3);

\draw[blue, thick] (1.5, 2.5) -- (2, 2.5);
\draw[blue, thick] (2, 2.5) .. controls +(right:3mm) and +(up:3mm) .. (2.5, 2);
\draw[blue, thick] (2.5, 2) .. controls +(down:3mm) and +(left:3mm) .. (3, 1.5);
\draw[blue, thick] (3, 1.5) .. controls +(right:3mm) and +(down:3mm) .. (3.5, 2);
\draw[blue, thick] (4, 2.5) .. controls +(left:3mm) and +(up:3mm) .. (3.5, 2);
\draw[blue, thick] (4, 2.5) -- (4.5, 2.5);

\draw[red, thick] (1.5, 1.5) -- (2, 1.5);
\draw[red, thick] (2, 1.5) .. controls +(right:3mm) and +(down:3mm) .. (2.5, 2);
\draw[red, thick] (2.5, 2) .. controls +(up:3mm) and +(left:3mm) .. (3, 2.5);
\draw[red, thick] (3, 2.5) .. controls +(right:3mm) and +(up:3mm) .. (3.5, 2);
\draw[red, thick] (3.5, 2) .. controls +(down:3mm) and +(left:3mm) .. (4, 1.5);
\draw[red, thick] (4, 1.5) -- (4.5, 1.5);

\end{tikzpicture}
\caption{Case 1.2 (blue) and case 2 (red).}
\label{case2}
\end{center}
\end{figure}
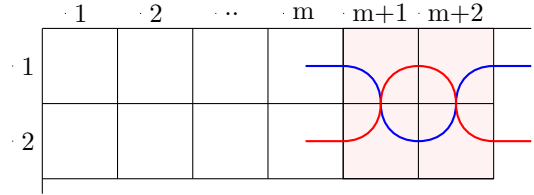

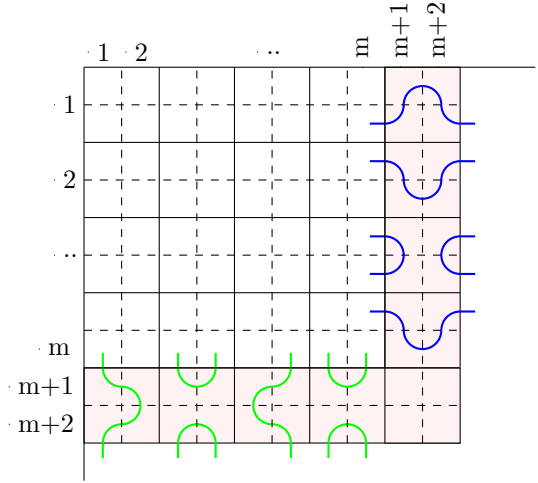
\begin{figure}[h]
\begin{center}
\begin{tikzpicture}

\draw (-1.95,3.2) -- (-1.95,3.2)   node[anchor=west]  {1};
\draw (-1.45,3.2) -- (-1.45,3.2)   node[anchor=west]  {2};
\draw (0.3,3.2) -- (0.3,3.2)   node[anchor=west]  {..};
\draw (1.7,3.0) -- (1.7,3.0)   node[rotate=90, anchor=west]  {m};
\draw (2.2,3.0) -- (2.2,3.0)   node[rotate=90, anchor=west]  {m+1};
\draw (2.7,3.0) -- (2.7,3.0)   node[rotate=90, anchor=west]  {m+2};
\draw (-2.4,2.5) -- (-2.4,2.5)   node[anchor=west]  {1};
\draw (-2.4,1.5) -- (-2.4,1.5)   node[anchor=west]  {2};
\draw (-2.4,0.5) -- (-2.4,0.5)   node[anchor=west]  {..};
\draw (-2.6,-0.75) -- (-2.6,-0.75)   node[anchor=west]  {m};
\draw (-3.0,-1.25) -- (-3.0,-1.25)   node[anchor=west]  {m+1};
\draw (-3.0,-1.75) -- (-3.0,-1.75)   node[anchor=west]  {m+2};

\filldraw[fill=red!5] (2,-2) rectangle (3,3);
\filldraw[fill=red!5] (-2,-2) rectangle (3,-1);


\draw[blue, thick] (1.8, 2.25) -- (2, 2.25);
\draw[blue, thick] (2, 2.25) arc (270:360:0.25);
\draw[blue, thick] (2.75, 2.5) arc (0:180:0.25);
\draw[blue, thick] (2.75, 2.5) arc (180:270:0.25);
\draw[blue, thick] (3, 2.25) -- (3.2, 2.25);

\draw[blue, thick] (1.8, 1.75) -- (2, 1.75);
\draw[blue, thick] (2.25, 1.5) arc (0:90:0.25);
\draw[blue, thick] (2.25, 1.5) arc (180:360:0.25);
\draw[blue, thick] (3, 1.75) arc (90:180:0.25);
\draw[blue, thick] (3, 1.75) -- (3.2, 1.75);


\draw[blue, thick] (1.8, -0.25) -- (2, -0.25);
\draw[blue, thick] (2.25, -0.5) arc (0:90:0.25);
\draw[blue, thick] (2.25, -0.5) arc (180:360:0.25);
\draw[blue, thick] (3, -0.25) arc (90:180:0.25);
\draw[blue, thick] (3, -0.25) -- (3.2, -0.25);


\draw[blue, thick] (1.8, 0.25) -- (2, 0.25);
\draw[blue, thick] (2, 0.25) arc (270:360:0.25);
\draw[blue, thick] (2.25, 0.5) arc (0:90:0.25);
\draw[blue, thick] (1.8, 0.75) -- (2, 0.75);
\draw[blue, thick] (3, 0.25) -- (3.2, 0.25);
\draw[blue, thick] (3, 0.75) arc (90:270:0.25);
\draw[blue, thick] (3, 0.75) -- (3.2, 0.75);

\draw[green, thick] (-1.75, -2.2) -- (-1.75, -2);
\draw[green, thick] (-1.5, -1.75) arc (90:180:0.25);
\draw[green, thick] (-1.5, -1.75) arc (270:360:0.25);
\draw[green, thick] (-1.25, -1.5) arc (0:90:0.25);
\draw[green, thick] (-1.75, -1) arc (180:270:0.25);
\draw[green, thick] (-1.75, -1) -- (-1.75, -0.8);

 

\draw[green, thick] (0.75, -2.2) -- (0.75, -2);
\draw[green, thick] (0.75, -2) arc (0:90:0.25);
\draw[green, thick] (0.5, -1.25) arc (90:270:0.25);
\draw[green, thick] (0.5, -1.25) arc (270:360:0.25);
\draw[green, thick] (0.75, -1) -- (0.75, -0.8);

\draw[green, thick] (-0.75, -2.2) -- (-0.75, -2);
\draw[green, thick] (-0.25, -2.2) -- (-0.25, -2);
\draw[green, thick] (-0.25, -2) arc (0:180:0.25);
\draw[green, thick] (-0.75, -1) -- (-0.75, -0.8);
\draw[green, thick] (-0.25, -1) -- (-0.25, -0.8);
\draw[green, thick] (-0.75, -1) arc (180:360:0.25);

\draw[green, thick] (1.25, -2.2) -- (1.25, -2);
\draw[green, thick] (1.75, -2.2) -- (1.75, -2);
\draw[green, thick] (1.75, -2) arc (0:180:0.25);
\draw[green, thick] (1.25, -1) -- (1.25, -0.8);
\draw[green, thick] (1.75, -1) -- (1.75, -0.8);
\draw[green, thick] (1.25, -1) arc (180:360:0.25);
\draw (-2,-1) -- (3,-1);
\draw (-2,0) -- (3,0);
\draw (-2,1) -- (3,1);
\draw (-2,2) -- (3,2);
\draw (-2,3) -- (4,3);
\draw[dashed] (-2,-1.5) -- (3,-1.5);
\draw[dashed] (-2,-0.5) -- (3,-0.5);
\draw[dashed] (-2,0.5) -- (3,0.5);
\draw[dashed] (-2,1.5) -- (3,1.5);
\draw[dashed] (-2,2.5) -- (3,2.5);
\draw (-2,-2.5) -- (-2,3);
\draw (-1,-2) -- (-1,3);
\draw (0,-2) -- (0,3);
\draw (1,-2) -- (1,3);
\draw (2,-2) -- (2,3);
\draw[dashed] (-1.5,-2) -- (-1.5,3);
\draw[dashed] (-0.5,-2) -- (-0.5,3);
\draw[dashed] (0.5,-2) -- (0.5,3);
\draw[dashed] (1.5,-2) -- (1.5,3);
\draw[dashed] (2.5,-2) -- (2.5,3);

\end{tikzpicture}
\caption{Example how the columns and rows $m+1$ and $m+2$ have to be crossed if no straight lines are allowed.}
\label{binned}
\end{center}
\end{figure}

\vspace{-5mm}
The rest of the proof is very similar to the proof of Lemma 6 with the only difference that we bin the chessboard in 2 by 2 squares. Like in Figure \ref{fig2}, where access to square $a_{1,m+1}$ was blocked and therefore the square below had to be passed on a horizontal path we now also cannot connect from the squares $a_{3,m+1}$ and $a_{3,m+2}$ to the squares above. We will have the same choices for the 4 squares $a_{3,m+1}$, $a_{3,m+2}$, $a_{4,m+1}$ and $a_{3,m+2}$ as we have for the 4 squares above which we just analysed. And mirrored, exactly as in the proof of Lemma 6, we will have no access from the squares $a_{m+1, 3}$ and $a_{m+1, 4}$ to the squares to the left if we do not allow straight lines in row $m+1$ and row $m+2$. Figure~\ref{binned} illustrates a typical situation where the analysed parts of columns and rows $m+1$ and $m+2$ are filled with possible paths.\\

\vspace{-5mm}

But what has to happen in the squares $a_{m+1,m+1}$, $a_{m+1,m+2}$, $a_{m+2,m+1}$ and $a_{m+2,m+2}$? It is the same situation as in the 2 x 2 field in the top left corner, since there is no access from left or from above. And this means, exactly as in the initialisation of the induction, that we need to have at least one straight line in either $a_{m+1,m+2}$ or $a_{m+2,m+1}$. And this is in contradiction to our assumption.\\

Therefore the $2(m+1) \times 2(m+1)$ corner must contain at least $m+1$ straight lines which completes the induction proof of Lemma 7.\\

\subsection{Final synthesis of the proof}

For $n$ being even, the minimum number of straight lines in a rook circuit on a $2n \times 2n$ chessboard follows directly from Lemma 7. We subdivide the chessboard in four quadrants of size $n \times n$ and apply the lemma in each quadrant. It follows immediately that the minimum number of straight lines is $4 \times \frac{n}{2} = 2n$.\\

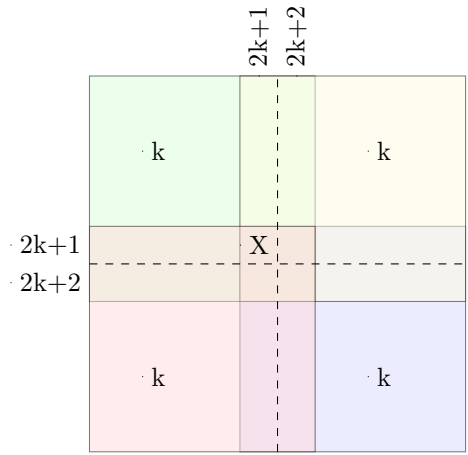
\begin{figure}[h]
\begin{center}
\begin{tikzpicture}

\filldraw[fill=blue!15, opacity=0.5] (-0.5,-2.5) rectangle (2.5,0.5);
\filldraw[fill=green!15, opacity=0.5] (-2.5,-0.5) rectangle (0.5,2.5);
\filldraw[fill=yellow!15, opacity=0.5] (-0.5,-0.5) rectangle (2.5,2.5);
\filldraw[fill=red!15, opacity=0.5] (-2.5,-2.5) rectangle (0.5,0.5);

\draw[dashed] (-2.5,0) -- (2.5,0);
\draw[dashed] (0,-2.5) -- (0,2.5);
\draw (-1.8,-1.5) -- (-1.8,-1.5)   node[anchor=west]  {k};
\draw (1.2,-1.5) -- (1.2,-1.5)   node[anchor=west]  {k};
\draw (-1.8,1.5) -- (-1.8,1.5)   node[anchor=west]  {k};
\draw (1.2,1.5) -- (1.2,1.5)   node[anchor=west]  {k};

\draw (-3.55,0.25) -- (-3.55,0.25)   node[anchor=west]  {2k+1};
\draw (-3.55,-0.25) -- (-3.55,-0.25)   node[anchor=west]  {2k+2};
\draw (-0.25,2.5) -- (-0.25,2.5)   node[rotate=90, anchor=west]  {2k+1};
\draw (0.25,2.5) -- (0.25,2.5)   node[rotate=90, anchor=west]  {2k+2};

\draw (-0.5,0.25) -- (-0.5,0.25)   node[anchor=west]  {X};

\end{tikzpicture}
\caption{A rook circuit on a $(4k+2) \times (4k+2)$ chessboard requires $4k + 4$ straight lines. X indicates a straight line outside the $2k \times 2k$ corners. Wherever the X is placed, three more straight lines have to be added to comply with the Lemmas 3 and 7.}
\label{quadrants}
\end{center}
\end{figure}

For $n$ being odd, we need to resort to Lemma 3, because we cannot subdivide the chessboard in four quadrants with a side length being even. But if we define $n = 2k + 1$ we can subdivide the chessboard in four overlapping $2(k+1) \times 2(k+1)$ corners. The four corners of size $2k \times 2k$ contain at least $k$ straights each. The situation is sketched in Figure \ref{quadrants}. According to Lemma 7, each $2(k+1) \times 2(k+1)$ corner must contain an additional straight compared to the $2k \times 2k$ corner. Since the four areas are overlapping, it is not necessarily required that we need one additional straight in each corner. But we need at least two, since the total number has to be even.\\

Let us try to add only two extra straights. There are two possibilities: either we will have two quadrants with $k+1$ straights and two with $k$ straights or we will have one quadrant with $k+2$ straights and three with $k$ straights. In the first case there is always a combination of two adjacent quadrants with a total of $2k+1$ straights. But the number of straights in two adjacent quadrants must be even because of Lemma 3. Hence this case can be excluded.\\

This leaves us with the case that we will have one quadrant with $k+2$ straights. Without loss of generality we assume that this is the green upper left quadrant. The blue lower right $2(k+1) \times 2(k+1)$ corner, diagonally opposed to the green quadrant, needs another straight. And therefore one of the straights in the green quadrant must be in the square $a_{2k+1,2k+1}$ (indicated with an "X" in Figure \ref{quadrants}) if we do not want to add extra straights.\\

But because of Lemma 3, we need an even number of straights both in column $2k+1$ as well as in row $2k+1$, which means we need at least another two straights to make the two numbers even. Hence also in this case, two additional straights are not sufficient.\\

Adding three additional straights violates the requirement that the number has to be even. Therefore we need to add at least four additional straights which means there will be at least $4k + 4$ straights in total.\\

In summary if $n$ is even the minimum number of straight lines is $2n$, if $n$ is odd the minimum number of straight lines is $2n + 2$.\\

q.e.d.\\

\subsection{Example solutions}

It is a nice exercise for the reader to find rook circuits with the minimum number of  straights. If $n$ is even, solutions are very easy to construct, because the number of columns and rows are multiples of 4 and you can start in one corner and spiral towards the centre. For instance, the rook meanders in columns $j$ with $j \bmod 4 = 3$ or 0 upwards, in columns $j$ with $j \bmod 4 = 1$ or 2 downwards, in rows $i$ with $i \bmod 4 = 3$ or 0 to the right and in rows $i$ with $i \bmod 4 = 1$ or 2 to the left. The rook arrives at half way in the central 4 x 4 square, which it can cross nicely in exactly the same pattern.\\

\begin{figure}[h!]
  \hspace{-15mm}
  \includegraphics[width=1.2\textwidth]{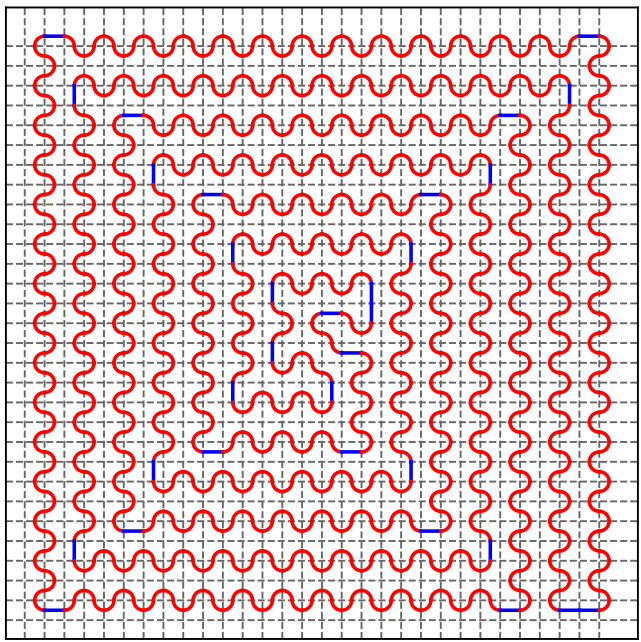}
  
  \vspace{-15mm}
  \caption{A rook circuit on a 30 x 30 chessboard with 32  straights (shown in blue)}
  \label{chessboard30x30}
\end{figure}

If $n$ is odd, the situation is slightly more complicated. Again the rook starts spiraling in from one corner, but when it arrives at the centre, a 6 x 6 square remains and the rook has to deviate from its strategy. Figure \ref{chessboard30x30} shows a solution how the central 6 x 6 square can be crossed with only 6 straight lines and thus limiting the circuit to the minimum straights.\\

\section{Minimum number of turns in a rook circuit on an $n \times m$ chessboard}

Let $k$ be the minimum number of turns in a rook circuit on an $n \times m$ chessboard with $n < m$, then 

$$
k = \left\{
    \begin{array}{ll}
        2n & \mbox{if } n \mbox{ is even} \\
        2m & \mbox{if } n \mbox{ is odd}
    \end{array}
\right.
$$

Proof:\\

Case 1: $n$ is even. We can split the chessboards in 4 corners of size $n/2 \times n/2$ and a middle strip of size $n \times (m-n)$. Each corner contains at least $n/2$ turns according to Lemma 6. Hence the total number of turns is at least $2n$.\\

Case 2: $n$ is odd. Without loss of generality we assume that we have $m$ rows and $n$ columns. A row cannot consist of $n$ straights, that would violate Lemma 4, because it would imply that we have an odd number of entries and exits in that row. So every row must contain at least two turns according to Lemma 3. Since we have $m$ rows and at least 2 turns in each row, the minimum number of turns is $2m$.\\

q.e.d.\\

The reader is invited to draw the circuits with the minimum number of turns. The circuit in Figure \ref{fig1} can be used as starting point. If $n$ is even, then we add $n-4$ rows containing to the right paths of this form: $\supset$ and stretch the horizontal lines to form $m$ squares. If $n$ is odd, then we add $m-4$ rows containing to the right paths of this form: $\supset$ and stretch the horizontal lines to form $n$ squares.

\section{Minimum number of straights in a rook circuit on an $n \times (n+1)$ chessboard}

Let $k$ be the minimum number of straights in a rook circuit on an $n \times (n+1)$ chessboard, then 

$$
k = \left\{
    \begin{array}{ll}
        n & \mbox{if } n \bmod 4 = 0 \\
        n + 1 & \mbox{if } n \bmod 4 = 1 \\
        n + 2 & \mbox{if } n \bmod 4 = 2 \\
        n + 1 & \mbox{if } n \bmod 4 = 3 
    \end{array}
\right.
$$

This corresponds to the OEIS sequence A201629 \cite{A201629}. 

\subsection{Proof for $n \bmod 4 = 0$} 

In this case the proof follows directly from Lemma 3. If $n = 4 k$ then we have without loss of generality $4k$ columns and $4m+1$ rows. In each column we need at least one straight in order to have an even number of turns in this column. Hence we will have at least $4k$ or $n$ straights.

\subsection{Proof for $n \bmod 4 = 1$} 

If $n = 4 k + 1$ then we have without loss of generality $4k + 2$ columns and $4k+1$ rows. Again from Lemma 3 it follows that we need in each column at least one straight in order to have an even number of turns in this column. Hence we will have at least $4k+2$ or $k+1$ straights. \\

\vspace{-4mm}
\subsection{Proof for $n \bmod 4 = 3$} 

If $n = 4 k - 1$ then we have without loss of generality $4k$ columns and $4k-1$ rows. Again from Lemma 3 it follows that we need in each column at least one straight in order to have an even number of turns in this column. Hence we will have at least $4k$ or $n+1$ straights. \\

\subsection{Proof for $n \bmod 4 = 2$} 

This is the only difficult case, because there will always be one column with more than one straight (assuming again an even number of columns and an odd number of rows).\\

\begin{figure}[h]
\begin{center}
\begin{tikzpicture}

\filldraw[fill=blue!15, opacity=0.5] (0.5,-3) rectangle (2.5,-1);
\filldraw[fill=green!15, opacity=0.5] (-2.5, 0.5) rectangle (-0.5,2.5);
\filldraw[fill=yellow!15, opacity=0.5] (0.5,0.5) rectangle (2.5,2.5);
\filldraw[fill=red!15, opacity=0.5] (-2.5,-3) rectangle (-0.5,-1);

\draw[dashed] (-2.5,0.5) -- (2.5,0.5);
\draw[dashed] (-2.5,0) -- (2.5,0);
\draw[dashed] (-2.5,-0.5) -- (2.5,-0.5);
\draw[dashed] (-2.5,-1) -- (2.5,-1);

\draw[dashed] (-0.5,-3) -- (-0.5,2.5);
\draw[dashed] (0,-3) -- (0,2.5);
\draw[dashed] (0.5,-3) -- (0.5,2.5);
\draw (-2.5, 0.5) -- (-2.5, -1);
\draw (2.5, 0.5) -- (2.5, -1);
\draw (-0.5, -3) -- (0.5, -3);
\draw (-0.5, 2.5) -- (0.5, 2.5);
\draw (-1.8,-2) -- (-1.8,-2)   node[anchor=west]  {kX};
\draw (1.2,-2) -- (1.2,-2)   node[anchor=west]  {kX};
\draw (-1.8,1.5) -- (-1.8,1.5)   node[anchor=west]  {kX};
\draw (1.2,1.5) -- (1.2,1.5)   node[anchor=west]  {kX};

\draw (-3.55,0.25) -- (-3.55,0.25)   node[anchor=west]  {2k+1};
\draw (-3.55,-0.25) -- (-3.55,-0.25)   node[anchor=west]  {2k+2};
\draw (-3.55,-0.75) -- (-3.55,-0.75)   node[anchor=west]  {2k+3};
\draw (-0.25,2.5) -- (-0.25,2.5)   node[rotate=90, anchor=west]  {2k+1};
\draw (0.25,2.5) -- (0.25,2.5)   node[rotate=90, anchor=west]  {2k+2};

\draw (-0.5,-0.25) -- (-0.5,-0.25)   node[anchor=west]  {X};
\draw (0,-0.25) -- (0,-0.25)   node[anchor=west]  {X};

\draw[red, thick] (-2.25, -1.2) -- (-2.25, -1);
\draw[blue, thick] (-1.75, -1) arc (0:90:0.25);
\draw[blue, thick] (-2.25, -0.5) arc (180:270:0.25);
\draw[blue, thick] (-2.0, -0.25) arc (90:180:0.25);
\draw[blue, thick] (-2.0, -0.25) arc (270:360:0.25);
\draw[blue, thick] (-1.75, 0) arc (0:90:0.25);
\draw[blue, thick] (-2.25, 0.5) arc (180:270:0.25);
\draw[blue, thick] (-2.25, 0.5) -- (-2.25, 0.7);

\draw[blue, thick] (-1.75, -1.2) -- (-1.75, -1);
\draw[red, thick] (-2, -0.75) arc (90:180:0.25);
\draw[red, thick] (-2, -0.75) arc (270:360:0.25);
\draw[red, thick] (-1.75, -0.5) arc (0:90:0.25);
\draw[red, thick] (-2, 0.25) arc (90:270:0.25);
\draw[red, thick] (-2, 0.25) arc (270:360:0.25);
\draw[red, thick] (-1.75, 0.5) -- (-1.75, 0.7);

\end{tikzpicture}
\caption{A rook circuit on a $(4k+2) \times (4k+3)$ chessboard requires $4k + 4$ straight lines. The X indicate where straight lines are required to be. For this chessboard two more straights will be required somewhere.}
\label{quadrantsn+1}
\end{center}
\end{figure}
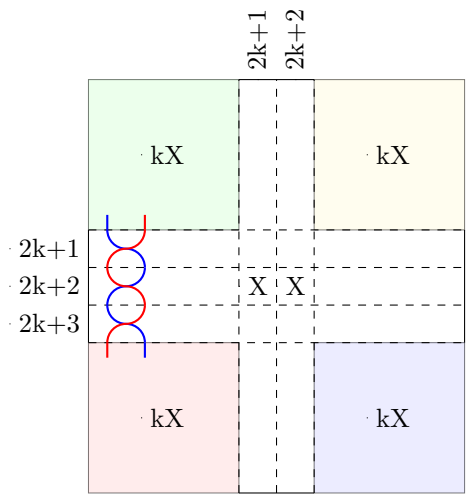

Let us assume for the moment we can find a rook circuit with only one straight per column, i.e.\ with $n$ straights in total. The only possibility to distribute the straights and complying with Lemmas 3 and 7 is shown in Figure~\ref{quadrantsn+1}. There are already $4k$ straights ($n=4k+2$) in the four $2k \times 2k$ corners according to Lemma 7. According to Lemma 3, there must be two more straights in columns $2k+1$ and $2k+2$, one in each column. Also from Lemma 3 it follows that they have to be in the same row (the number of straights per row must be even) and from Lemma 7 it follows that they have to be in the middle row (otherwise there are two $(2k+2) \times (2k+2)$ corners that will not have at least $k+1$ straights). Hence they are in $a_{2k+2,2k+1}$ and $a_{2k+2,2k+2}$, indicated with an $X$. \\

\vspace{-7mm}

And now we can use the same reasoning that we used for the proof of Lemma 7. Figure \ref{binned} showed how the 4 squares $a_{m+1,1}$, $a_{m+1,2}$, $a_{m+2,1}$ and $a_{m+2,2}$ had to be crossed if straights must be avoided. The same logic can be applied to show that the 6 squares $a_{2k+1,1}$, $a_{2k+1,2}$, $a_{2k+2,1}$, $a_{2k+2,2}$, $a_{2k+3,1}$ and $a_{2k+3,2}$ can only be crossed on two ways (on the red or blue path shown in Figure \ref{quadrantsn+1}). The same pattern is now enforced in column 3 and 4 and this proliferates until columns $2k-1$ and $2k$. The same reasoning is applied starting from the right most columns towards the left.\\

This brings us to the situation where the 6 central squares in columns $2k+1$ and $2k+2$ cannot be connected with squares to the left or right. Hence the two $X$ must be vertical straights. But what does that mean for the four squares $a_{2k+1,2k+1}$, $a_{2k+1,2k+2}$, $a_{2k+3,2k+1}$ and $a_{2k+3,2k+2}$ right above and below the two $X$? Since they cannot be straights and they cannot go to the outside, they must all go to the inside and the six central squares must be connected in one small separate circuit. This of course violates the requirement of one single rook circuit and hence the assumption that the rook circuit can be constructed with only $n$ straights collapses.

\subsection{Example solutions}

\begin{figure}[h!]
  \hspace{-1mm}
  \includegraphics[width=1.05\textwidth]{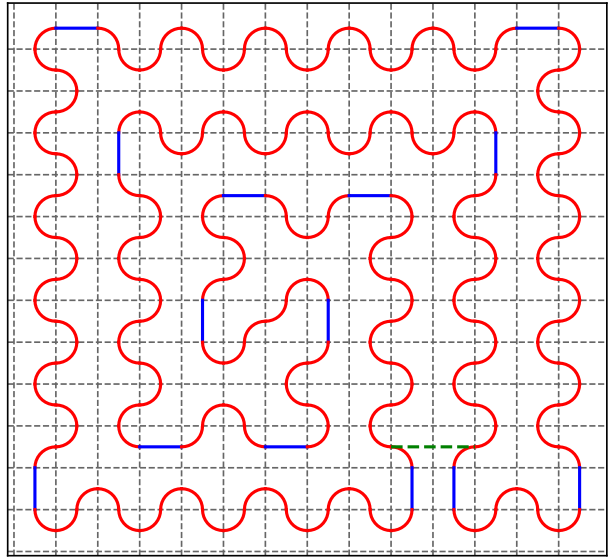}
  
  \vspace{-10mm}
  \caption{A rook circuit on a 13 x 14 chessboard with 14  straights (shown in blue). If you replace two turns by the green dashed line you can see the solution for the 9 x 10 chessboard.}
  \label{chessboard13x14}
\end{figure}

It is a nice exercise for the reader to find rook circuits with the minimum number of straights. For most cases solutions are pretty easy to find, but it took me two days to find a solution for the 9 x 10 chessboard. Certain solutions for an $n \times m$ board, can be easily extended to arrive at a solution for the $(n+4) \times (m+4)$ board. Figure \ref{chessboard13x14} shows a solution for a 13 x 14 chessboard. It is based on the solution of the 9 x 10 chessboard. Only $a_{11,10}$ and $a_{11,11}$ need to be exchanged with horizontal straights (as indicated with the green dashed line) and the solution of the 9 x 10 board becomes visible.\\

  \vspace{-7mm}

Also this solution for the 13 x 14 board can easily be extended for a solution of the 17 x 18 board by only modifying the 8 squares in the bottom right corner and adding one slalom spiral around the whole circuit. Figure \ref{chessboard30x30} nicely shows how more and more spirals can be added to a previous solution.\\ 

\section{Minimum number of straights in a rook circuit on an $n \times m$ chessboard}

In the case that $n \bmod 4 = 0$, it is easy to draw circuits with only $n$ straights. And according to Lemma 7 it cannot be done with less straights. Also in case that $n \bmod 4 = 1$, it is easy to draw circuits with $m$ straights which is the minimum following the logic of Lemma 3. The conjectured results for $n \bmod 4 = 2$ and 3 are also listed in Table \ref{table:1}, but the cases with a gray shading are not straightforward to prove. This is left for future work.

\begin{table}[h!]
\centering
\begin{tabular}[h]{|c|c|c|c|c|}
\hline
\diagbox{n mod 4}{m mod 4} & 0 & 1 & 2 & 3 \\
\hline
0 & n & n & n & n \\
1 & m & - & m & - \\
2 & \cellcolor{Gray}m & \cellcolor{Gray} m-1 & \cellcolor{Gray} m+2 & \cellcolor{Gray} m+1 \\
3 & m & - & \cellcolor{Gray} m+2 & - \\
\hline 
\end{tabular}
\caption{Minimum number of straights in a rook circuit on an $n \times m$ chessboard. The gray shaded results still remain to be proved.}
\label{table:1}
\end{table}

\bibliography{references} 
\bibliographystyle{plain} 

\end{document}